\documentclass{amsart}

\usepackage{amssymb,amsthm,amsmath}   

 \newcommand{\bee}{\begin{equation}}
 \newcommand{\eee}{\end{equation}}

\newcommand{\diam}{\mbox{\rm diam}}
\newcommand{\be}{\begin{eqnarray}}
\newcommand{\ee}{\end{eqnarray}}
\newcommand{\supp}{\mbox{\rm supp}}

\newcommand{\const}{\mbox{\rm const}}
\newcommand{\Span}{\mbox{\rm Span}}

\newcommand{\eps}{{\mbox{$\epsilon$}}}

\newcommand{\Th}{{\theta}}
\newcommand{\gam}{{\gamma}}

\newcommand{\R}{{\mathbb R}}
\newcommand{\Q}{{\mathbb Q}}
\newcommand{\Z}{{\mathbb Z}}
\def\Gbb{{\mathbb G}}
\newcommand{\C}{{\mathbb C}}
\newcommand{\Nat}{{\mathbb N}}

\newcommand{\Ak}{{\mathcal A}}

\newcommand{\Pk}{{\mathcal P}}

\newcommand{\Ck}{{\mathcal C}}

\newcommand{\Kk}{{\mathcal K}}

\newcommand{\Vk}{{\mathcal V}}

\newcommand{\Sk}{{\mathcal S}}
\newcommand{\Tk}{{\mathcal T}}

\newcommand{\Uk}{{\mathcal U}}

\newcommand{\dist}{\mbox{\rm dist}}

\newcommand{\lam}{\lambda}

\newcommand{\om}{\omega}

\newcommand{\rank}{{\rm rank}}

 \newtheorem{theorem}{Theorem}[section]
 \newtheorem{lemma}[theorem]{Lemma}
 
 \newtheorem{cor}[theorem]{Corollary}

\theoremstyle{definition}
\newtheorem{defi}[theorem]{Definition}

\theoremstyle{remark}

\numberwithin{equation}{section}

\begin{document}

\title[Eigenfunctions for tiling substitutions]{Eigenfunctions for substitution tiling systems}

\author[B. Solomyak]{Boris Solomyak$^{\rm 1}$}

\address{Boris Solomyak, Box 354350, Department of Mathematics,
University of Washington, Seattle WA 98195}
\email{solomyak@math.washington.edu}

\thanks{$^1$Rosi and Max Varon Visiting Professor at the
Weizmann Institute of Science in the Autumn of 2005.
This research was supported in
part by NSF grant DMS 0355187.}

\begin{abstract}
We prove that for the uniquely ergodic $\R^d$-action associated with a
primitive substitution tiling of finite local complexity,
every measurable eigenfunction coincides with a  continuous function 
almost everywhere.
Thus, topological weak-mixing is equivalent to measure-theoretic
weak-mixing for such actions.
If the expansion map for the substitution is a pure dilation by $\theta>1$
and the substitution has a fixed point, then failure of  
weak-mixing is equivalent to $\theta$ being a Pisot number.
\end{abstract}

\maketitle

\section{Introduction}

In this article we consider self-affine (substitution) tilings of $\R^d$ and 
associated dynamical systems. These tilings are of translationally finite
local complexity,
that is, they have a finite number of tiles and patches of a given size,
up to translation.
A self-affine tiling has the property that if we inflate it by
a certain expanding linear map,
then the original tiling may be obtained by subdividing
the inflated tiles according to a prescribed rule.
Self-affine
tilings have been extensively studied;
they arise, in particular, in connection
with Markov partitions of toral automorphisms and as models for quasicrystals.
See the next section for some historical comments.

We focus on the dynamical system, which may be
defined (under appropriate assumptions)
as the translation action by $\R^d$ on the orbit closure of  a
given tiling. 
Assuming primitivity of the tiling substitution,
we obtain a uniquely ergodic action. The spectrum (more precisely,
the spectral measure of the associated unitary group) is a fundamental
invariant of the dynamical system, and it is closely related to the
diffraction spectrum studied by physicists (see \cite{Dw,LMS1}).
In particular, the point spectrum (the set of eigenvalues) 
corresponds to the Bragg peaks (``sharp bright spots'') in the diffraction
picture associated with the tiling.

The new result (announced in \cite{Sol2}) proved
here is that every measurable eigenfunction can be chosen to be continuous.
Thus, the topological
and ergodic-theoretic point spectra are the same (such systems are
sometimes called homogeneous \cite{Robi1}). This extends a 
result of Host \cite{Host} on symbolic substitution $\Z$-actions.
The condition that all eigenfunctions are continuous has been used recently
in the work on mathematical quasicrystals, see \cite{BLM,L}.

The characterization and existence of eigenvalues have a link with
Number Theory, namely with Pisot (or PV) numbers. We present a proof of
the following statement: if the expanding map is a pure dilation
by $\theta>1$, then the dynamical system has non-trivial eigenvalues
(that is, it is not weak-mixing) if and only if $\theta$ is Pisot.
A similar result was obtained 
in \cite{GK} in the context of diffraction
spectrum, by different methods.

\medskip

{\bf Acknowledgment.} This is an expanded version of 
the talk given at the International Conference on Probability and Number
Theory held in Kanazawa in June 2005. I am grateful to  the
organizers for their warm hospitality. 
Thanks also to Rick Kenyon for useful discussions.

\section{Historical Remarks}

Here we make some historical remarks, in order to put our work in a larger
context. They are necessarily incomplete, and we apologize for any
inadvertent omissions. We refer to the survey by E. A. Robinson, Jr.
\cite{Robi} for more details.

The story begins in mathematical logic
\cite{Wan61,Brg66} with the discovery of 
{\em aperiodic prototile sets}, i.e.\ 
sets of prototiles which can tile the plane, but only in a such a way that
the resulting tiling does not have any translational symmetries.
One of the most interesting aperiodic prototile sets was
constructed in the early 1970's by R. Penrose \cite{Penrose}.
Penrose tilings have an additional feature that they can be generated using
a tiling substitution. See \cite{GS} for a detailed description of Penrose
tilings and many other aperiodic sets.

In 1984 physicists discovered what came to be known as {\em quasicrystals}
\cite{SBGC}.
These are metallic alloys, which, like a crystal, have a sharp
X-ray diffraction pattern, but unlike a crystal, have an aperiodic atomic
structure. Aperiodicity was inferred from the ``forbidden'' 5-fold symmetry
of the diffraction pattern. 
It turned out that Penrose tilings have similar features,
so they became a focus of many investigations, both by
physicists and by mathematicians. Other tilings have been studied from this
point of view as well.
See \cite{Senechal,baake} for an introduction to
the mathematics of quasicrystals and further references.

The third source of our subject is symbolic dynamics (and closely
related word combinatorics) and ergodic theory.
Prototiles may be considered as a kind of geometric symbols and entire
tilings as infinite (multi-dimensional for $d>1$) words.
It is useful to consider not just an individual tiling,
but a space of tilings, together with the translation action.
A new feature is that this is an action by $\R^d$, rather than $\Z^d$, and
topologically,
the tiling space is not a Cantor set, but a solenoid. Tiling dynamical
systems were first introduced by D. Rudolph \cite{Rud}, and further
investigated
by C. Radin and M. Wolff \cite{RW} and E. A. Robinson, Jr. \cite{Robi0},
see also the book \cite{radbook}. As already mentioned in the Introduction,
the dynamical spectrum is closely related to the diffraction spectrum
(see \cite{Dw,LMS1}), and this connection has been used to obtain new results
about diffraction, see e.g.\ \cite{Schlot}.

Tilings with an inflation symmetry, of which the Penrose tiling is
an example, may be viewed as generalizations of substitution sequences.
Substitutions (also called morphisms) have been studied
in many fields, including ergodic theory and
dynamical systems, see e.g.\ \cite{Pytheas} and references therein.
The spectral theory of substitution dynamical systems was
developed by M. Queffelec \cite{Queff}, B. Host \cite{Host} and other 
authors.
Independently, {\em adic transformations} were introduced by Vershik \cite{Ver}
as a general model for measure-preserving systems. Their theory
was developed
by A. M. Vershik, A. Livshits, and others (see \cite{Liv0,VerLiv} and references
therein), which, in the stationary case,
was linked to the theory of substitutions \cite{Liv1,HZ}.
The importance of Pisot numbers for the spectral properties of substitution
and adic systems was suggested by A. Livshits in the mid 1980's
[personal communication]. B. Host discussed it in 
\cite[(6.2)]{Host}. A complete algebraic characterization of eigenvalues for
substitution systems, using the notion of Pisot families,
was obtained by S. Ferenczi, C. Mauduit and A. Nogueira \cite{FMN}.
A geometric realization of a certain Pisot substitution as a ``domain exchange''
on a torus, which also yielded a self-similar tiling of the plane with
fractal boundaries, was discovered by
G. Rauzy \cite{Rau}, which inspired a lot of new research.
An axiomatic framework was introduced by
W. Thurston \cite{Thur}, who studied self-similar tilings of the plane and 
their expansion constants. This theory was further developed by R. Kenyon
\cite{Ken,KenSST}. Self-affine tilings have been used for constructions of
Markov partitions, see e.g.\ \cite{Bedford,Ken,prag,IO,KenVer}.
There are also strong links with numeration systems and beta-expansions,
see \cite{BeSie} and references therein.

Although it had precursors, such as \cite{RW,BeRad},
it seems that \cite{soltil} was the first systematic attempt to extend
substitution dynamics to the tiling setting. It was continued in \cite{Sol2}
which extended B. Moss\'e's results \cite{Mosse} on recognizability.
The present paper is a further study in this direction.

Finally, we should add that there are many other interesting developments
in the subject which go far beyond the scope of this paper, 
among them the study of $C^*$-algebras arising from 
substitution tilings \cite{AP}, deformations of tiling spaces \cite{ClaSad},
conjugacies for tiling dynamical systems \cite{HRS},
tiling dynamical systems as $\Gbb$-solenoids and laminated spaces \cite{BeGa},
etc.


\section{Preliminaries and Statement of Results}

We begin with tiling preliminaries, following \cite{LMS2,Robi}.
We emphasize that our tilings are {\em translationally
finite}, thus excluding the pinwheel tiling \cite{radin} and its relatives.

\subsection{Tilings.}
Fix a set of types (or colors) labeled by $\{1,\ldots,m\}$.
A {\em tile} in $\R^d$ is defined as a pair $T = (A,i)$ where
$A  = \supp(T)$ (the support of $T$) is a compact set in $\R^d$ which is the
closure of its interior, and $i = \ell(T)\in \{1,\ldots,m\}$ is the type of $T$.
(The tiles are not assumed to
be homeomorphic to the ball or even connected.)
A {\em tiling} of $\R^d$ is a set $\Tk$ of tiles such that $\R^d = \bigcup
\{\supp(T):\ T\in \Tk\}$ and distinct tiles (or rather, their supports)
have disjoint interiors.

A patch $P$ is a finite set of tiles with disjoint interiors.
The {\em support of a patch} $P$ is defined by
$\supp (P) = \bigcup\{\supp(T):\ T\in P\}$. The {\em diameter of a patch}
$P$ is $\diam(P)=\diam(\supp(P))$. The {\em translate} of a tile $T=(A,i)$
by a vector $g\in \R^d$ is $T+g = (A+g, i)$. The translate
of a patch $P$ is $P+g = \{T+g:\ T\in P\}$. We say that two patches
$P_1,P_2$ are {\em translationally equivalent} if $P_2 = P_1+g$ for some
$g\in \R^d$. Finite subsets of $\Tk$ are called $\Tk$-patches.

\begin{defi} \label{def-flc}
A tiling $\Tk$ has (translational) {\em finite local complexity} (FLC) if
for any $R>0$ there are finitely many $\Tk$-patches of diameter less
than $R$ up to translation equivalence.
\end{defi}

\begin{defi} \label{def-rep}
A tiling $\Tk$ is called {\em repetitive} if for any patch $P\subset \Tk$
there is $R>0$ such that for any $x\in \R^d$ there is a
$\Tk$-patch $P'$ such that $\supp(P')\subset B_R(x)$ and $P'$ is a
translate of $P$.
\end{defi}

\subsection{Tile-substitutions, self-affine tilings.}
We study {\em perfect}
(geometric) substitutions, in which a tile is ``blown up'' by an expanding
linear map and then subdivided. Other possibilities, where the substitution
is combinatorial, and/or there is no perfect geometric subdivision, have also
been considered, see e.g.\ \cite{Priebe1,Priebe2,PriSol,ClaSad}.

A linear map $\phi : \R^d \rightarrow \R^d$ is {\em expansive}
if all its eigenvalues lie outside the unit circle.

\begin{defi}\label{def-subst}
Let $\Ak = \{T_1,\ldots,T_m\}$ be a finite set of tiles in $\R^d$
such that $T_i=(A_i,i)$; we will call them {\em prototiles}.
Denote by $\Pk_{\Ak}$ the set of
patches made of tiles each of which is a translate of one of $T_i$'s.
A map $\omega: \Ak \to \Pk_{\Ak}$ is called a {\em tile-substitution} with
expansion $\phi$ if 
\begin{equation} \label{def-sub}
\supp(\om(T_j)) = \phi A_j \ \ \  \mbox{for} \  j\le m.
\end{equation}

In plain language, every expanded prototile $\phi T_j$ can be decomposed into
a union of tiles (which are all translates of the prototiles) with disjoint
interiors.
\end{defi}

The substitution $\om$ is extended to all translates of prototiles by
$\om(x+T_j)= \phi x + \om(T_j)$, and to patches by
$\om(P)=\cup\{\om(T):\ T\in P\}$. This is well-defined due to (\ref{def-sub}).
The substitution $\om$ also acts on the space of tilings whose tiles are
translates of those in $\Ak$.

To the substitution $\om$ we associate its $m \times m$
substitution matrix ${\sf S}$, with ${\sf S}_{ij}$ 
being the number of tiles of type $i$
in the patch $\om(T_j)$.
The substitution $\om$ is called {\em primitive}
if the substitution matrix is primitive, that is, if there exists $k\in \Nat$
such that ${\sf S}^k$ has only positive entries.
We say that $\Tk$ is a fixed point of a substitution if $\om(\Tk) = \Tk$.

\begin{defi} \label{def-tsp}
Given a primitive tile-substitution $\om$, let
$X_\om$ be the set of all tilings whose
every patch is a translate of a subpatch of
$\om^n(T_j)$ for some $j\le m$ and $n\in \Nat$. (Of course, one can use
a specific $j$ by primitivity.) The set $X_\om$ is called the
{\em tiling space} corresponding to the substitution.
\end{defi}

\begin{defi} \label{def-saf}
A repetitive FLC fixed point of a primitive tile-substitution 
is called a {\em self-affine tiling}.
It is called {\em self-similar} if the expansion map is a similitude, that is,
$|\phi(x)| = \Th |x|$ for all $x\in \R^d$, with some $\Th>1$.
\end{defi}

It is often convenient to
work with self-affine tilings; doing this is not a serious restriction,
since every primitive tile-substitution has a periodic point in the tiling
space and replacing $\om$ by $\om^n$ does not change the
tiling space.

We say that a tile-substitution $\om$ has FLC
if for any $R>0$ there are finitely many subpatches of $\om^n(T_j)$
for all $j\le m$, $n\in \Nat$, of diameter less
than $R$, up to translation. This obviously implies that all tilings in
$X_\om$ have FLC, and is equivalent to it if the tile-substitution is
primitive.

\medskip

{\bf Remark.}
A primitive substitution tiling
space is not necessarily of finite local complexity,
see \cite{Ken1,Danzer,FraRob}. Thus we have
to assume FLC explicitly.

\begin{lemma}\cite[Prop.\,1.2]{prag} \label{lem-prag}
Let $\om$ be a primitive tile-substitution
of finite local complexity.
Then every tiling $\Sk\in X_\om$ is repetitive.
\end{lemma}

\subsection{Tiling topology and tiling dynamical system} 
We use a tiling metric on $X_\om$, which is based on a
simple idea: two tilings are close if after a small translation they agree
on a large ball around the origin. There is more than one way
to make this precise.
We say that two tilings $\Tk_1,\Tk_2$ {\em agree} on a set $K \subset \R^d$ if
$$
\supp(\Tk_1\cap \Tk_2) \supset K.
$$
For $\Tk_1,\Tk_2 \in X_\om$ let
$$
\widetilde{d}(\Tk_1,\Tk_2) := \inf\{r \in (0,2^{-1/2}):
\ \exists\,g,\ \|g\| \le r\ 
\mbox{such that}
$$ 
$$
\Tk_1-g \mbox{ \ agrees with $\Tk_2$ on $B_{1/r}(0)$} \}.
$$
Then
$$
d(\Tk_1,\Tk_2) = \min\{2^{-1/2},\widetilde{d}(\Tk_1,\Tk_2)\}.
$$

\begin{theorem}\cite{Rud} {\em (see also \cite{Robi}).} 
$(X_\om,d)$ is a complete metric space.
It is compact, whenever the space has finite local complexity. 
The action of $\R^d$ by translations on
$X_\om$, given by  $g(\Sk)= \Sk-g$, is continuous.
\end{theorem}

This continuous translation action $(X_\om,\R^d)$ is called the
(topological) tiling dynamical system associated with the tile-substitution.

\begin{theorem} \label{th-min}
If $\om$ is a primitive tiling substitution with FLC, then
the dynamical system $(X_\om,\R^d)$ is minimal, that is, for every
$\Sk\in X_\om$, the orbit $\{\Sk-g:\ g\in \R^d\}$ is dense in $X_\om$.
\end{theorem}

This follows from Lemma~\ref{lem-prag} and Gottschalk's Theorem
\cite{Gott}, see \cite[Sec.\,5]{Robi} for details.

\begin{defi} 
A vector $\alpha =(\alpha_1,\ldots,\alpha_d) \in \R^d$ is said to be
an eigenvalue for the continuous $\R^d$-action if there exists an eigenfunction
$f\in C(X_\om)$, that is, $\ f\not\equiv 0$ and for all $g\in \R^d$ and
all $\Sk\in X_\om$,
\begin{equation} \label{def-eig1}
f(\Sk-g) = e^{2 \pi i \langle g, \alpha\rangle} f(\Sk).
\end{equation}
Here $\langle \cdot,\cdot \rangle$ denotes the standard scalar product
in $\R^d$.
\end{defi}

Note that this ``eigenvalue'' is actually a vector. In physics it might be
called a ``wave vector.'' More generally, for an action of a locally compact
Abelian group $G$, the eigenvalues are elements of the dual group
$\widehat{G}$. For a single transformation (translation by a vector $g$),
the eigenvalue is the more familiar $e^{2 \pi i \langle g, \alpha\rangle}$,
a point on the unit circle.

\subsection{Measurable dynamics.}
A topological dynamical system is said to be uniquely ergodic if it has
a unique invariant Borel probability measure. 

\begin{theorem} \label{unerg}
If $\om$ is a primitive tiling substitution with FLC, then
the dynamical system $(X_\om,\R^d)$ is uniquely ergodic.
\end{theorem}

This result has appeared in the literature in several
slightly different versions.
We refer to \cite{LMS2,Robi} for the proof.

Let $\mu$ be the unique invariant measure from Theorem~\ref{unerg}.
The measure-preserving tiling dynamical system is denoted by 
$(X_\om,\R^d,\mu)$.

\begin{sloppypar}
\begin{defi} \label{def-eig2}
A vector $\alpha \in \R^d$ is an eigenvalue for the
measure-preserving system $(X_\om,\R^d,\mu)$ if
there exists an eigenfunction
$f\in L^2(X_\om,\mu)$, that is, $f$ is not the zero function in $L^2$
and for all $g\in \R^d$,
the equation (\ref{def-eig1}) holds for $\mu$-a.e.\ $\Sk\in X_\om$.
\end{defi}
\end{sloppypar}

By ergodicity, all the eigenvalues are simple and the eigenfunctions 
have a constant modulus a.e., see \cite{Wal}.
To distinguish between the measure-theoretic and topological settings, we
can speak about measurable and continuous eigenfunctions. 

\begin{theorem} \label{th-main}
If $\om$ is a primitive tiling substitution with FLC, then every
measurable eigenfunction for the system $(X_\om,\R^d,\mu)$ coincides with a
continuous function $\mu$-a.e.
\end{theorem}

This extends the result of Host \cite{Host} on $\Z$-actions
associated to primitive one-dimensional symbolic substitutions.

\medskip

{\bf Remark.} Continuous and measurable eigenfunctions for
linearly recurrent Cantor systems were
recently investigated in \cite{CDHM,BDM}. In the latter paper necessary
and sufficient conditions for being an eigenvalue are established, and it
is proved that not every measurable eigenfunction is a.e.\ continuous for such
systems.

\medskip

A dynamical system is said to be {\em weak-mixing} if it has no
non-constant eigenfunctions. This notion is considered both in the
topological and the measure-theoretic category.
As a consequence of Theorem~\ref{th-main}, for 
our systems measure-theoretic weak-mixing is equivalent to 
topological weak-mixing.

Theorem~\ref{th-main} was announced in \cite{Sol2}. 
In the case of aperiodic substitution spaces, the result is immediate from
\cite[Th.\,5.1]{soltil} and \cite[Th.\,1.1]{Sol2}, so that here we only need to
deal with the case when some periods are present.

\subsection{Characterization of eigenvalues.}
Let $\om$ be a primitive tiling substitution of finite local complexity and
let $\Sk\in X_\om$.
Consider the set of translation vectors between tiles of the same type:
\begin{equation} \label{def-xi}
\Xi: = \{x\in \R^d:\ \exists\,T,T'\in \Sk,\ T'= T+x\}.
\end{equation}
It is clear that $\Xi$ does not depend on the tiling $\Sk$. The vectors
$x\in \Xi$ are sometimes called {\em return vectors}; they are tiling analogs
of return words in word combinatorics. We also need the
group of translation symmetries
$$
\Kk = \{x\in \R^d:\ \Sk-x = \Sk\},
$$
which does not depend on $\Sk$ either. The tiling space is said to be
{\em aperiodic} if $\Kk = \{0\}$, {\em sub-periodic} if 
$0 < \rank (\Kk) < d$, and {\em periodic} if $\rank (\Kk) = d$.

\begin{theorem} \label{th-char} 
Let $\om$ be a primitive tiling substitution of finite local complexity
with expansion map $\phi$, which has a fixed point (a self-affine tiling).
Then the following are equivalent for $\alpha \in \R^d$:

{\bf (i)} $\alpha$ is an eigenvalue for the topological dynamical system
$(X_\om,\R^d)$;

{\bf (ii)} $\alpha$ is an eigenvalue for the measure-preserving system
$(X_\om,\R^d)$;

{\bf (iii)} $\alpha$ satisfies the following two conditions:
\begin{equation} \label{eig1}
\lim_{n\to\infty} e^{2\pi i \langle \phi^n z, \alpha \rangle} =1
\ \ \mbox{for all}\ z\in \Xi,
\end{equation}
and
\begin{equation} \label{eig2}
e^{2\pi i \langle g, \alpha \rangle} =1
\ \ \mbox{for all}\ g\in \Kk.
\end{equation}
\end{theorem}
 
This theorem is also a generalization of the corresponding result from
\cite{Host}.
Theorem~\ref{th-main} is immediate from Theorem~\ref{th-char}, since 
eigenvalues are simple for an ergodic system and hence normalized 
measurable and continuous eigenfunctions for the same eigenvalue
must coincide a.e.

\medskip

Theorem~\ref{th-char} does not address the question when the eigenvalues are
present, so additional analysis is needed. This turns out to be closely
related to Number Theory, more precisely, to {\em Pisot numbers} (also
called PV-numbers) and their
generalizations. Recall that an algebraic integer $\Th>1$ is a Pisot number
if all its Galois conjugates $\theta'$
(other roots of the minimal polynomial) satisfy $|\theta'|<1$. 

\begin{theorem} \label{th-Pisot}
Let $\theta>1$, and let $\phi(x) = \theta x$ on
$\R^d$ be the expansion map. Let $\om$ be a primitive tile-substitution of
finite local complexity with
expansion $\phi$, admitting a fixed point. 
Then the associated measure-preserving system
is not weak-mixing if and only if $\theta$ is a Pisot number.
\end{theorem}

A similar result was obtained by G\"ahler and Klitzing
\cite{GK}, where the diffraction
spectrum was considered, from a different point of view. We 
present a proof using our methods. Both \cite{GK} and our
proof rely on a result of R. Kenyon, which says
that under the assumptions of the theorem, $\Xi \subset \Z[\theta]b_1 + \cdots
+ \Z[\theta]b_d$, for some basis $\{b_1,\ldots,b_d\}$ of $\R^d$.
We include a proof of the latter as well, since it is not easy to
extract from the literature.

\medskip

{\bf Remarks.}
1. We do not have a general theorem like Theorem~\ref{th-Pisot},
for an arbitrary expansion map
$\phi$, but some partial results
 which involve complex Pisot numbers and
Pisot families, may be found in \cite{soltil,Robi}.

2. It is interesting to compare the simple criterion for weak mixing
in Theorem~\ref{th-Pisot} with the case of 
substitution $\Z$-actions, for which a more complicated
characterization was obtained
by S. Ferenczi, C. Mauduit and A. Nogueira \cite{FMN}
The reason, roughly, is that here we only 
consider ``geometric'' tiling substitutions,
so the algebraic conjugates of the expansion constant do not enter into the
picture.

3. In this paper we do not address the difficult question whether the
spectrum is {\em pure discrete}, see \cite{soltil,BD,LMS2}, as well as the
more recent \cite{BK,L} and references therein.


\section{Continuous eigenfunctions}

In this section we prove Theorem~\ref{th-char}. To deduce
Theorem~\ref{th-main}, we note that for a primitive tile-substitution $\om$
there exist $\Tk\in X_\om$ and $k\in \Nat$ such
that $\om^k(\Tk) = \Tk$ (see \cite[Th.\,5.10]{Robi}).
Since
replacing $\om$ by $\om^k$ does not 
change the tiling space ($\phi$ should be replaced by $\phi^k$),
the result would follow.

\medskip

The implication (i) $\Rightarrow$ (ii) in Theorem~\ref{th-char} is obvious. 

\medskip

{\em Proof of} (ii) $\Rightarrow$ (iii).
The necessity of (\ref{eig1}) is proved in \cite[Th.\,4.3]{soltil}.
We should note that in \cite{soltil} it was assumed that the expansion map
$\phi$ is diagonalizable, but the proof of \cite[Th.\,4.3]{soltil} works 
for any expansion map.
The straightforward details are left to the interested reader.

Next we prove the necessity of (\ref{eig2}). Let $g\in \Kk$. Then
$\Sk -g = \Sk$ for every $\Sk \in X_\om$. If $f_\alpha$ is a measurable
eigenfunction corresponding to $\alpha \in \R^d$, then
$$
f_\alpha(\Sk) = f_\alpha(\Sk-g) = e^{2 \pi i \langle g, \alpha\rangle} f(\Sk)
$$
for $\mu$-a.e.\ $\Sk\in X_\om$. It follows that 
$e^{2 \pi i \langle g, \alpha\rangle}=1$, as desired.

\medskip

{\em Proof of} (iii) $\Rightarrow$ (i).
Let $\Tk$ be the self-affine tiling whose existence we assumed. Suppose that
(\ref{eig1}), (\ref{eig2}) hold and define
\begin{equation} \label{Eig}
f_\alpha(\Tk - x) = e^{2\pi i \langle x,\alpha \rangle}\ \ \mbox{for}\ 
x\in \R^d.
\end{equation}
The orbit $\{\Tk-x:\ x\in \R^d\}$ 
is dense in $X_\om$ by minimality. If we show that $f_\alpha$ is
uniformly continuous on this orbit, then we can extend $f_\alpha$ to $X_\om$,
and this extension will satisfy the eigenvalue equation (\ref{def-eig1}) by
continuity. We will need several lemmas; the first one will be useful in the
next section as well.

\begin{lemma} \cite{Ken,Thur}
\label{lem-expa1} Suppose that there is a primitive tiling
substitution of finite local complexity with expansion $\phi$. Then
all the eigenvalues of $\phi$ are algebraic integers.
\end{lemma}

{\em Proof.} We provide a short proof for completeness.
Consider $\langle \Xi \rangle$, the subgroup of
$\R^d$ generated by $\Xi$.
It is a free finitely generated Abelian group by
FLC. Thus, we can find a set of free generators $v_1,\ldots,v_\ell$ for
$\langle \Xi \rangle$. 
Consider the $d\times \ell$ matrix $V =
[v_1 \ldots v_\ell]$.
By the definition of substitution tilings, $\phi \Xi \subset \Xi$, hence
$\phi$ acts on $\langle \Xi \rangle$, and so there is an integer
$\ell\times \ell$
matrix $M$ such that 
\begin{equation} \label{eq-defM}
\phi V = VM.
\end{equation}
Note that $\ell\ge d$ and $\rank(V)=d$ since $\Xi$ spans $\R^d$.
If $\lam$ is an eigenvalue of $\phi^T$ with the eigenvector $x$, then
$\lam$ is also an eigenvalue of $M^T$ with the eigenvector $V^Tx$.
(The superscript $T$ denotes the transpose of a matrix.)
Since $M$ is an integer matrix, it follows that $\lam$ is an algebraic
integer. \qed

\begin{lemma} \label{lem-exp} If (\ref{eig1}) holds for some $z\in \Xi$,
then the convergence in (\ref{eig1}) is exponential; in fact, there exist
$\rho \in (0,1)$, depending only on the expansion $\phi$, and $C_z >0$ such 
that
\begin{equation} \label{eq-expo}
\left| e^{2\pi i \langle \phi^n z, \alpha \rangle} - 1\right| < C_z \rho^n
\ \ \mbox{for }\ n\in \Nat.
\end{equation}
\end{lemma}

Below we will show that (\ref{eq-expo}) holds with a constant
$C$ independent of $z$, but this lemma is the first step.

\medskip

{\em Proof.} This can be deduced from \cite[L.\,2]{Maud},
but we sketch a (well-known) direct proof for the reader's convenience.

We continue the argument in the proof of the previous lemma.
Let $z\in \Xi$. By the definition of
free generators, there is a unique vector $a(z) \in \Z^\ell$ such that
$z = V a(z)$. Then we have from (\ref{eq-defM}):
$$
\langle \phi^n z,\alpha \rangle = \langle \phi^n Va(z),\alpha \rangle
= \langle VM^na(z),\alpha \rangle = 
\langle M^na(z),V^T\alpha \rangle:= \zeta_n.
$$
It follows from the Caley-Hamilton Theorem
that the sequence $\zeta_n$ satisfies a recurrence relation
with integer coefficients:
\begin{equation} \label{eq-recur}
\zeta_{n+\ell} + a_1 \zeta_{n+\ell-1} + \cdots + a_n \zeta_n = 0,
\end{equation}
where $t^\ell + a_1 t^{\ell-1} 
+ \cdots + a_\ell$ is the characteristic polynomial $p$ of
the integer matrix $M$.
Let 
$$
\zeta_n = \langle \phi^n z,\alpha \rangle := K_n + \eps_n,
$$
where $K_n$ is the nearest integer to $\zeta_n$. By (\ref{eig1}), we have
$\eps_n \to 0$.
Since the sequence $\zeta_n$ satisfies (\ref{eq-recur}) 
and $K_n$ are integers,
we conclude that $K_n$ satisfy the same recurrence relation for $n$
sufficiently large. Therefore, also $\eps_n$ satisfy (\ref{eq-recur}),
hence they can be expressed in terms of
the zeros of $p$,
for $n$ sufficiently large. More precisely, there exist
$\alpha_j \in \C$ and polynomials $q_j \in \Z[x]$ such that
\begin{equation} \label{eq-expa}
\eps_n = \sum_{j=1}^{\ell_1}
\alpha_j q_j(n) \theta_j^n \mbox{\ \  for}  \ n \ge n_0,
\end{equation}
where $\theta_j,\ j=1,\ldots,\ell_1$, are the distinct zeros of $p$.
Since $\eps_n\to 0$, it is not hard to see
that all $\theta_j$, which occur in (\ref{eq-expa}) with nonzero 
coefficients, must satisfy $|\theta_j|< 1$. Then
$$
|\eps_n| \le \const \cdot \rho^n,
$$

\begin{sloppypar}
\noindent where $\max\{|\theta_j|:\ |\theta_j| < 1\} < \rho < 1$, and
we conclude that
$\left| e^{2\pi i \langle \phi^n z, \alpha \rangle} - 1\right| =
\left| e^{2\pi i\eps_n} - 1\right|$ satisfies
(\ref{eq-expo}) for appropriate $C_z$. \qed
\end{sloppypar}

\medskip

Next we need a lemma which is analogous to \cite[Th.\,1.5]{prag} and
\cite[Lem.\,6.5]{soltil}. We do not include a complete proof, but deduce it
from \cite{LS}.

\begin{lemma} \label{lem-exp2}
Let $\om$ be a primitive tile-substitution of finite local complexity with
expansion $\phi$, which has a fixed point $\Tk$,
and let $\Xi$ be the set of translation vectors between
tiles of the same type in $X_\om$. Then there exists $k\in \Nat$
and a finite set
$U$ in $\R^d$, with $\phi^k U \subset \Xi$, such that for any $z\in \Xi$
there exist $N\in \Nat$ and $u(j), w(j) \in U,\ 0 \le j \le N$, so that
$$
z = \sum_{j=0}^N \phi^{kj}(u(j) + w(j)).
$$
\end{lemma} 

{\em Proof.} We can find $k\in \Nat$ so that ${\sf S}^k$ is strictly positive,
where ${\sf S}$ is the substitution matrix of $\om$.
Then $\Tk$ is also a fixed point of $\om^k$, a tile-substitution with
expansion $\phi^k$ and a strictly positive substitution matrix. Now we can
apply \cite[Lem.\,4.5]{LS} to obtain the desired result. \qed

\begin{cor} \label{cor-unif} If (\ref{eig1}) holds, then the convergence
is uniform in $z \in \Xi$, that is,
$$
\lim_{n\to\infty} \sup_{z\in \Xi} 
\left| e^{2\pi i \langle \phi^n z, \alpha \rangle} - 1\right| = 0.
$$
\end{cor}

{\em Proof.} Let
$$
C := \max\{C_u:\ u\in U\},
$$
where $C_u$ is from Lemma~\ref{lem-exp} and $U$ is from Lemma~\ref{lem-exp2}.
Let $z\in \Xi$ and consider the expansion from Lemma~\ref{lem-exp2}.
Then we have
\begin{eqnarray*}
\left|e^{2\pi i \langle \phi^n z,\alpha \rangle} - 1 \right| & = &
\Bigl|\exp\Bigl(2\pi i \Bigl\langle \sum_{j=0}^N \phi^{n+kj} (u(j) + w(j)),
\alpha \Bigr\rangle \Bigr) -1 \Bigr| \\
& \le & \sum_{j=0}^N \left| \exp(2\pi i \langle \phi^{n+k(j-1)}
\phi^k u(j),\alpha
 \rangle ) -1 \right| \\
& + & \sum_{j=0}^N \left| \exp(2\pi i \langle \phi^{n+k(j-1)}\phi^k w(j),\alpha
 \rangle ) -1 \right| \\
& \le & 2C \sum_{j=0}^N \rho^{n+k(j-1)}  < 2C\rho^{-k} (1-\rho^k)^{-1}\rho^n,
\end{eqnarray*}
which implies the desired statement. Here we used that
$|e^{i(a+b)} - 1| \le |e^{ia}-1| + |e^{ib}-1|$ for real $a,b$ 
and (\ref{eq-expo}).
\qed

\medskip

Extending the work of Moss\'e \cite{Mosse} on primitive aperiodic substitution
$\Z$-actions,
we proved in \cite{Sol2} that
if $X_\om$ is aperiodic, then the substitution map $\om$ is bijective on
$X_\om$. In the general case, we proved in \cite[Th.\,1.2]{Sol2} that
if $\om(\Sk_1) = \om(\Sk_2)$, then $\Sk_2 = \Sk_1 - \phi^{-1}g$ for some $g\in
\Kk$. The next lemma is a quantitative version of this
statement, which follows from it easily.

\begin{lemma} \label{lem-invert} Let $\om$ be a primitive tile-substitution
of finite local complexity. 
Then for every $\eps>0$ there is $\delta>0$ such that
for any $\Sk,{\Sk}'\in X_\om$ with $d(\Sk,{\Sk}') < \delta$
there exist $\Uk,{\Uk}'\in X_\om$ such that
$\om(\Uk) = \Sk,\ \om({\Uk}') = {\Sk}'$, and $d(\Uk,{\Uk}') < \eps$.
\end{lemma}

{\em Proof.}
Recall that $\Kk$ is the set of periods for $X_\om$, so that
$\Sk-h = \Sk$ for all $h\in \Kk$. Consider the 
following equivalence relation on $X_\om$:
$$
\Sk_1 \sim \Sk_2 \ \Longleftrightarrow\ \Sk_2 = \Sk_1-\phi^{-1}g
\ \ \mbox{for some}\ g \in \Kk.
$$
Observe that $\Kk$ is a discrete subgroup of $\R^d$ and $\phi\Kk$ is
a subgroup of finite index in $\Kk$. Thus, all the equivalence classes
are finite. Denote the equivalence class of $\Sk$ by $[\Sk]$,
and the set of equivalence classes by $\widehat{X}_\om$.
Consider the induced metric on $\widehat{X}_\om$:
$$
\widehat{d}([ \Sk_1 ]\,, [ \Sk_2 ]):=
\min\{d(\Sk_1',\Sk_2'):\ \Sk_1' \sim \Sk_1,\ \Sk_2'\sim\Sk_2\}.
$$
It is readily seen that $(\widehat{X}_\om, \widehat{d})$ is compact.
Consider the map 
$\widehat{\om}:\,\widehat{X}_\om \to X_\om$ given by
$$
\widehat{\om}([ \Vk ]) = \om(\Vk)\ \ \mbox{for}\ 
\Vk\in X_\om.
$$
This is well-defined, since for $\Vk'\sim \Vk$ we have $\Vk' = \Vk-g$ for
some $g\in \phi^{-1}\Kk$, hence $\om(\Vk') = \om(\Vk) - \phi g = \om(\Vk)$.
Since $\om$ is a continuous surjection onto $X_\om$, we have that
$\widehat{\om}$ is a continuous surjection from 
$\widehat{X}_\om$ onto $X_\om$.
We claim that $\widehat{\om}$ is 1-to-1. Indeed, if 
$\widehat{\om}([ \Sk_1 ])=
\widehat{\om}([ \Sk_2 ])$,
then $\om(\Sk_1) = \om(\Sk_2)$, and $\Sk_1 \sim \Sk_2$ by 
\cite[Th.\,1.2]{soltil}. It follows that $\widehat{\om}^{-1}$ is 
uniformly continuous, which is precisely the desired statement.
\qed

\medskip

{\em Conclusion of the proof of} (iii) $\Rightarrow$ (i) {\em in 
Theorem~\ref{th-char}.}
Recall that we need to show the uniform continuity of $f_\alpha$, given by
(\ref{Eig}), on the orbit $\{\Tk-x:\ x\in \R^d\}$. We have $\om(\Tk) = \Tk$,
so for any $h\in \R^d$, by \cite[Th.\,1.2]{soltil},
$$
\om^{-1}(\Tk-h) = \{\Tk - \phi^{-1}h - \phi^{-1}g:\ g \in \Kk\}.
$$
Applying Lemma~\ref{lem-invert}, we obtain that for any $\eps>0$ there exists
$\delta>0$ such that
$$
\mbox{if \ }\ d(\Tk- x, \Tk - y) < \delta\ \ \mbox{then}
$$
\begin{equation} \label{duk1}
\exists\,g,g'\in \Kk,\
d(\Tk - \phi^{-1}x - \phi^{-1} g, \Tk - \phi^{-1}y - \phi^{-1} g') < \eps.
\end{equation}

Fix $\eta\in (0,1)$. 
By Corollary~\ref{cor-unif}, we can choose $n\in \Nat$ such that
$|e^{2\pi i \langle \phi^n z,\alpha \rangle} - 1| < \eta/2$ for all
$z\in \Xi$. Applying (\ref{duk1}) $n$ times, we can find $\delta>0$ such that
if $d(\Tk- x, \Tk - y) < \delta$, then there exist $g_j,g_j'\in \Kk$, for
$1 \le j \le n$, with
$$
d\Bigl(\Tk - \phi^{-n}x -\sum_{j=1}^n \phi^{-n + j -1} g_j,
\Tk - \phi^{-n}y - \sum_{j=1}^n \phi^{-n + j -1} g_j'\Bigr) < 
\eps
$$
where
$$
\eps:=\frac{\eta}{4\pi\|\phi^n\|\|\alpha\|}\,.
$$
By the definition of the metric $d$, this means that there exists
$h\in \R^d$, with $\|h\| \le \eps$, such that 
$
\Tk - \phi^{-n}x -\sum_{j=1}^n \phi^{-n + j -1} g_j-h$
agrees with
$
\Tk - \phi^{-n}y - \sum_{j=1}^n \phi^{-n + j -1} g_j'$ on
$B_{1/\eps}(0)$.
Agreement on any neighborhood of the origin implies that the tilings share
the tiles containing the origin. Thus, by  the definition of the set $\Xi$,
$$
\phi^{-n} (x-y) + \sum_{j=1}^n \phi^{-n + j -1}(g_j - g_j') - h \in \Xi.
$$
Therefore, there exists $z\in \Xi$ such that
$$
x-y =  \phi^n z + \phi^n h + w,\ \  \ \mbox{where}\ w:= 
\sum_{j=1}^n \phi^{j -1}(g'_j - g_j) \in \Kk
$$
(here we use that $\phi \Kk\subset \Kk$).
Finally,
\begin{eqnarray*}
|f_\alpha(x) - f_\alpha(y)| & = &
\left|e^{2\pi i \langle x-y,\alpha \rangle} - 1 \right| \\
& = & \left|e^{2\pi i (\langle \phi^nz,\alpha \rangle+
\langle \phi^n h,\alpha \rangle + \langle w, \alpha \rangle)} - 1 \right| \\
& \le & \left|e^{2\pi i \langle \phi^nz,\alpha \rangle}-1\right| +
\left|e^{2\pi i\langle \phi^n h,\alpha \rangle} -1\right| \\
& \le & \eta/2 + 2\pi|\langle \phi^n h,\alpha \rangle| \\
& \le & \eta/2 + 2\pi\|\phi^n\|\|\alpha\|\eps = \eta,
\end{eqnarray*}
and we are done. We used the condition (\ref{eig2}) to get rid of the term 
with $w$.
 \qed


\section{Pisot substitutions}

{\em Proof of necessity in Theorem~\ref{th-Pisot}.}
We need to prove that if the dynamical system $(X_\om,\R^d)$ has a
non-constant eigenfunction, then $\theta$ is a Pisot number. Recall
that here
we assume the expansion map to be a pure dilation: $\phi(x) = \theta x$.

Let $\alpha \ne 0$ be an eigenvalue. The set of translation vectors
$\Xi$ between tiles of the same type spans $\R^d$, hence we can find
$z\in \Xi$ such that $\langle z, \alpha \rangle \ne 0$. By 
Theorem~\ref{th-char}, the distance from $\theta^n \langle z, \alpha \rangle$
to the nearest integer tends to zero, as $n\to \infty$. We know that
$\theta$ is algebraic (see Lemma~\ref{lem-expa1}), hence $\theta$ is a Pisot
number by the classical result of Pisot (see e.g. \cite{Salem}). \qed

\medskip

{\em Proof of sufficiency in Theorem~\ref{th-Pisot}.}
We need to show that if $\theta$ is Pisot, then there are non-zero eigenvalues
for the dynamical system. The proof relies on the following result.

\begin{theorem}[Kenyon] \label{th-ken} 
Let $\Tk$ be a self-similar tiling with expansion map $\phi(x) = \Th x$ for
some $\Th>1$ and let $\Xi$ be the set of return vectors of the tiling $\Tk$,
defined in (\ref{def-xi}).
Then there exists a basis $\{b_1,\ldots,b_d\}$ of $\R^d$
such that 
\begin{equation} \label{eq-ken}
\Xi \subset \Z[\theta]b_1 + \cdots + \Z[\theta]b_d.
\end{equation}
\end{theorem}

First we finish the proof of sufficiency.
Suppose that there are non-trivial periods, that is, $\Kk \ne \{0\}$.
Since $\theta \Kk \subset \Kk$ and $\Kk$ is a discrete subgroup of $\R^d$,
we obtain that $\theta \in \Nat$. Then it follows from (\ref{eq-ken}) that
the group generated by $\Xi$, which we denoted by $\langle \Xi \rangle$, is
discrete. It is a lattice in $\R^d$, since it spans $\R^d$. It is clear that
the all points of the dual lattice $\langle \Xi \rangle'$ satisfy
both (\ref{eig1}) and (\ref{eig2}) (using that $\Kk \subset \Xi$),
hence there are non-trivial eigenvalues.

Now suppose that the tiling is aperiodic, that is, $\Kk= \{0\}$.
Let $\{b_1^*,\ldots,b_d^*\}$ be the dual basis for  $\{b_1,\ldots,b_d\}$,
that is, $\langle b_i, b_j^* \rangle = \delta_{ij}$. We claim that the set
$$
\Z[\Th^{-1}]b_1^*  + \cdots + \Z[\Th^{-1}]b_d^*
$$

\begin{sloppypar}
\noindent is contained in the group of eigenvalues. Indeed, suppose
$\alpha = \sum_{j=1}^d b_j^* p_j(\Th^{-1})$ for some polynomials $p_j\in
\Z[x]$. Let
$z\in \Xi$. By (\ref{eq-ken}), we can write
$
z = \sum_{j=1}^d b_j q_j(\theta)
$
for some polynomials $q_j\in \Z[x]$. 
Then
\end{sloppypar}
$$
\langle \phi^n z, \alpha \rangle = \theta^n \sum_{j=1}^d q_j(\theta)
p_j(\theta^{-1}) = \theta^{n-k} P(\Th)
$$ 
for some $k\in \Nat$ and $P\in \Z[x]$. We have $\dist(\theta^{n-k} P(\Th),
\Z) \to 0$, as $n\to \infty$ (see \cite{Salem}), so (\ref{eig1}) is
satisfied and $\alpha$ is an eigenvalue by Theorem~\ref{th-char}. \qed

\medskip

{\em Proof of Theorem~\ref{th-ken}.}
Our proof is based on \cite{Thur,Ken2} and a 
personal communication from Rick Kenyon; however, we do not need quasiconformal
maps as in \cite{Ken2} which is concerned with more general tilings.

Instead of the set $\Xi$, it is more convenient to work with {\em control
points}, see \cite{Thur,Ken,prag}.

\begin{defi} \cite{prag}
Let $\Tk$ be a fixed point of a primitive substitution with expansive map
$\phi$. For each $\Tk$-tile $T$, fix a tile $\gamma T$ in the patch
$\omega (T)$;
choose $\gamma T$ with the same relative position for all tiles of the same
type. This defines a map $\gamma : \Tk \to \Tk$ called the
{\em tile map}. Then define the {\em control point} for a tile $T \in \Tk$ by
\[ c(T) = \bigcap_{n=0}^{\infty} \phi^{-n}(\gamma^n T).\]
Note that the control points are not uniquely defined; they depend on the
choice of $\gam$.
\end{defi}

Let $\Ck=
\Ck(\Tk) = \{c(T):\ T\in \Tk\}$ be the set of control points for all tiles.
The control points have the following properties:
\begin{itemize}
\item[(a)] $T' = T + c(T') - c(T)$, for any tiles $T, T'$ of the same type;
\item[(b)] $\phi(c(T)) = c(\gamma T)$, for $T \in \Tk$.
\end{itemize}
Therefore,
\begin{itemize}
\item[(c)] $\Xi \subset \Ck-\Ck$.
\item[(d)] $\phi(\Ck) \subset \Ck$.
\end{itemize}

The definition above works for a general expansion $\phi$, but
now we again assume that
$\phi(x) = \theta x$.
Observe that it is enough to prove the inclusion
\begin{equation} \label{eq-ken2}
\Ck \subset \Q(\theta)e_1 + \cdots +  \Q(\theta)e_d
\end{equation}
for some basis $\{e_1,\ldots, e_d\}$. Indeed, the Abelian group
$\langle \Ck\rangle$
is finitely generated. Let $\{w_1,\ldots,w_N\}$ be a set of free generators.
By (\ref{eq-ken2}), $w_j = \sum_{j=1}^d e_j 
\frac{p_j^{(i)}(\Th)}{q_j^{(i)}(\Th)}$, for
$i\le N$, for some polynomials $p_j^{(i)}, q_j^{(i)}\in \Z[x]$.
Then we obtain (\ref{eq-ken}), with
$b_j = e_j 
\left(\prod_{j=1}^d\prod_{i=1}^N q_j^{(i)}(\Th)\right)^{-1}$, in view of
the property (c) above.

\medskip

Now pick any set $\{e_1,\ldots,e_d\}\subset \Ck$ which
spans $\R^d$.
Consider the vector space 
$$
\widetilde{\Ck}:= \Span_{\Q(\Th)}\Ck
$$
over the field $\Q(\Th)$.
We want to show that $\{e_1,\ldots,e_d\}$ is a basis for 
$\widetilde{\Ck}$. Let $\pi$ be any linear projection from $\widetilde{\Ck}$
onto $\Span_{\Q(\Th)}\{e_1,\ldots,e_d\}$. Since the vector space is over
$\Q(\Th)$, we obviously have
\begin{equation} \label{eq-comm}
\pi(\Th^n x) = \Th^n \pi (x)\ \ \ 
\mbox{for}\ x\in \widetilde{\Ck},\ n\in \Z.
\end{equation}
Thus we have a map
\begin{equation} \label{eq-deff}
f(x) = \pi(x) \ \ \ \mbox{for}\ x\in \Ck_\infty,\ \  
\mbox{where}\ \ \Ck_\infty:= \bigcup_{n=0}^\infty \theta^{-n} \Ck.
\end{equation}
This is consistent in view of (\ref{eq-comm}).
Now $f$ is defined on a dense subset of $\R^d$; let us show that
it is uniformly Lipschitz on this set, then $f$ can be extended to $\R^d$
by continuity.
Fix any norm $\|\cdot\|$ on $\R^d$.

\begin{lemma} \label{lem-lip} 
There exists $L_1>0$ such that
\begin{equation} \label{eq-lip1}
\|\pi(\xi)-\pi(\xi')\| \le L_1 \|\xi - \xi'\|,\ \ \forall\ \xi,\xi'\in
\Ck.
\end{equation}
\end{lemma}

{\em Proof.} This argument is from \cite{Thur}.
One can move ``quasi-efficiently''
between control points by moving from ``neighbor to neighbor.'' More
precisely, let
$W$ be the set of vectors in $\R^d$ obtained as the
difference vectors between control points of
neighboring $\Tk$-tiles.
Note  that $W$ is finite by FLC.
There is a constant
$C_1=C_1(\Tk)$ such that $\forall\ \xi,\xi'\in\Ck$, there exist $p\in\Nat$ and
$\xi_1: = \xi,\ \xi_2,\ldots,\xi_{p-1}\in \Ck,\ \xi_p:=\xi'$
such that $\xi_{i+1}-\xi_i \in W$
for $i=1,\ldots,p-1$, and
$$
\sum_{i=1}^{p-1}\|\xi_{i+1}-\xi_i\|\le C_1\cdot \|\xi-\xi'\|.
$$
(This is an exercise; see \cite{Lag} for a detailed proof.)
Let
$$
C_2:= \max\{\|\pi(w)\|/\|w\|:\ w\in W\}.
$$
Now we can estimate:
\begin{eqnarray*}
\|\pi(\xi)-\pi(\xi')\| = \|\pi(\xi-\xi')\| & = & \left\|\sum_{i=1}^{p-1}
\pi(\xi_{i+1}-\xi_i)\right\|  \\
& \le & \sum_{i=1}^{p-1}\|\pi(\xi_{i+1}-\xi_i)\|
\\ & \le & C_2 \sum_{i=1}^{p-1}\|\xi_{i+1}-\xi_i\| \\
& \le & C_1 C_2 \|\xi-\xi'\|.
\end{eqnarray*}
\qed

\medskip

In view of (\ref{eq-comm}), the last lemma
implies that $f$ is Lipschitz on $\Ck_\infty$ with the uniform Lipschitz
constant $L_1$. We extend $f$ to a Lipschitz function on $\R^d$
by continuity; it satisfies
\begin{equation} \label{eq-comm2}
f(\Th x) = \Th f(x)\ \ \ \mbox{for all} \ x\in \R^d.
\end{equation}
Note that the extension $f$ need not coincide with $\pi$ on 
$\widetilde{\Ck}$; we started with $\pi$ not on the entire $\Q(\Th)$-vector
space, but only on the control points and their preimages under the
expansion.

\begin{lemma} \label{lem-loc}
The function $f$ depends only on the tile type in $\Tk$, up to an
additive constant:
if $T, T+x \in \Tk$ and $\xi \in \supp(T)$, then
\begin{equation} \label{eq-loc}
f(\xi+x) = f(\xi) + \pi(x).
\end{equation}
\end{lemma}


{\em Proof.}
Observe that $x\in \Xi \subset \Ck - \Ck \subset \widetilde{\Ck}$, so
$\pi(x)$ is defined.
It is enough to check (\ref{eq-loc}) on a dense set. Suppose
$\xi = \Th^{-k} c(S) \in \supp(T)$ for some $S\in \om^k(T)$. Then
$S + \Th^k x \in \om^k(T+x)\subset \Tk$ and we have
\begin{eqnarray*}
f(\xi+x) & = & f(\Th^{-k} c(S) + x) \\
& = & f (\Th^{-k} c(S + \Th^kx)) \\
& = & \Th^{-k} \pi(c(S + \Th^kx)) \\
& = & \Th^{-k} \pi(c(S)) + \Th^{-k} \pi (\Th^k x) \\
& = & f(\xi) + \pi(x),
\end{eqnarray*}
as desired. \qed

\medskip

{\em Conclusion of the proof of Theorem~\ref{th-Pisot}.}
We mimic the argument of Thurston \cite{Thur} but provide more details.

The function $f:\ \R^d\to \R^d$ is Lipschitz,
hence it is differentiable almost everywhere. Let $x$ be a point where
the total derivative $H=Df(x)$ exists. 
Then
$$
f(x+u) = f(x) + Hu + \Psi(u) \ \ \ \mbox{for all}\ \ u\in \R^d,
$$
where 
\begin{equation} \label{little-o}
\|\Psi(u)\|/\|u\| \to 0,\ \ \mbox{as}\ u\to 0.
\end{equation}
Multiplying by $\theta^k$, using (\ref{eq-comm2}) and substituting $v=\Th^ku$,
we obtain
$$
f(\Th^k x + v) = f(\Th^k x) + Hv + \Th^k\Psi(\Th^{-k}v)\ \ \ \mbox{for all}\ \
v\in \R^d.
$$
For a set $A\subset \R^d$ denote by $[A]^{\Tk}$ the $\Tk$-patch of $\Tk$-tiles
whose supports intersect $A$. By repetitivity, there exists $R>0$ such that
$B_R(0)$ contains a translate of the patch $[B_1(\Th^kx)]^{\Tk}$ for all
$k\in \Nat$. This implies, in view of Lemma~\ref{lem-loc}, that there exist
$x_k \in B_R(0)$, for $k\ge 1$, such that
$$
f(x_k+v) = f(x_k) + Hv + \Th^k\Psi(\Th^{-k}v)\ \ \ \mbox{for all}\ 
\ v\in B_1(0).
$$
There exists a limit point $x'$ of the sequence $\{x_k\}$; then by
continuity of $f$ and (\ref{little-o}),
$$
f(x'+ v) = f(x') + Hv\ \ \ \mbox{for all}\ \ v\in B_1(0).
$$
Thus, $f$ is flat on some neighborhood. Applying (\ref{eq-comm2}) again,
we obtain that $f$ is flat on an arbitrarily large neighborhood.
By repetitivity, a translate of $[B_1(0)]^{\Tk}$ occurs in every sufficiently
large neighborhood, therefore, by Lemma~\ref{lem-loc}, the function $f$ is
flat on $B_1(0)$. Using (\ref{eq-comm2}) for the last time, we conclude that
$f$ is flat everywhere, and since $f(0)=0$ (see (\ref{eq-comm}) and 
(\ref{eq-deff})), $f$ is linear.
But $f(e_j) = \pi(e_j)=e_j$ for the set $\{e_j\}_{j\le d}$
which we chose in the beginning of the proof, hence $f$ is the identity map,
$\pi(\xi) = \xi$ for all $\xi\in \Ck$, and we conclude that 
(\ref{eq-ken2}) holds, as desired. \qed

\end{document}